\documentclass[graybox]{svmult}
\usepackage{amsfonts, amsmath, amssymb}
\usepackage{mathtools}
\usepackage{mathrsfs}
\usepackage{mathptmx}
\usepackage{booktabs}
\usepackage{protosem}
\usepackage{caption}
\usepackage{epigraph}
\usepackage{enumitem}
\usepackage{subcaption}
\usepackage{tikz}
\usepackage{float}
\usepackage[numbers,sort&compress]{natbib}

\usepackage{lipsum}

\usepackage{todonotes}

\usepackage{hyperref}

\newcommand{\R}{\mathbb{R}}

\newcommand{\N}{\mathbb{N}}
\newcommand{\Z}{\mathbb{Z}}

\newtheorem{experiment}{\sc Experiment}
{\bf}{\it}

\newtheorem{defi}{Definition}

\usepackage{type1cm}        % activate if the above 3 fonts are
\usepackage{makeidx}         % allows index generation
\usepackage{graphicx}        % standard LaTeX graphics tool
\usepackage{multicol}        % used for the two-column index
\usepackage[bottom]{footmisc}% places footnotes at page bottom

\usepackage{newtxtext} 
\usepackage[varvw]{newtxmath}       % selects Times Roman as basic font
\makeindex             % used for the subject index
\begin{document}

\title*
{
% Hybrid Surrogate Models: Accurate Solutions of Non-Smooth Partial Differential Equations by Overcoming the Gibbs Phenomenon
Hybrid Surrogate Models: Circumventing Gibbs Phenomenon for Partial Differential Equations with Finite Shock-Type  Discontinuities
}

\titlerunning{Hybrid Surrogate Models} 
% for an abbreviated version of
\author{Juan-Esteban Suarez Cardona\orcidID{0000-0002-1990-2413} and\\
Shashank Reddy\orcidID{0009-0006-7037-6312} and \\
Michael Hecht\orcidID{0000-0001-9214-8253}}

\institute{Juan-Esteban Suarez Cardona \at Center for Advanced Systems Understanding (CASUS), G\"{o}rlitz, Germany, \email{j.suarez-cardona@hzdr.de}
\and Shashank Reddy \at Institute of Information Technology, Bengalore, India\email{Shashank.Reddy@iiitb.ac.in} \and Michael Hecht \at Center for Advanced Systems Understanding (CASUS), G\"{o}rlitz, Germany \at University of Wroclaw, Mathematical Institute,  Wroclaw, Poland \email{m.hecht@hzdr.de}}

%\institute{\email{,m.hecht@hzdr.de}
%\at 
%\at $^1$
%\at Helmholtz-Zentrum Dresden-Rossendorf e.V. (HZDR), Dresden, Germany$^2$
%\at University of Wroclaw, Mathematical Institute,  Wroclaw, Poland$^3$ 
%\at Technical University  Dresden,
%Faculty of Mathematics, Dresden, Germany$^4$
%\at , $^5$
%}

% \institute{Juan Esteban Suarez Cardona$^{1,2}$, 
% \email{j.suarez-cardona@hzdr.de}

% \and Shashank Reddy \at Institute of Information Technology, Bengalore, India 

% \and Michael Hecht$^{1,2,3}$, \email{m.hecht@hzdr.de} 
% \at Center for Advanced Systems Understanding (CASUS), G\"{o}rlitz, Germany$^1$
% \at 
% Helmholtz-Zentrum Dresden-Rossendorf e.V. (HZDR), Dresden, Germany$^2$ 
% \at Mathematical Institute, University of Wroclaw, Wroclaw, Poland$^3$
\maketitle
\vspace{-6em}\abstract{
We introduce the concept of \emph{Hybrid Surrogate Models} (HSMs) -- combining multivariate polynomials with Heavyside functions --  as approximates of functions with finitely many jump discontinuities.  We exploit the HSMs for formulating a variational optimization approach, solving non-regular partial differential equations (PDEs) with non-continuous shock-type solutions. 
The HSM technique simultaneously obtains a parametrization of the position and the height of the shocks as well as the solution of the PDE. We show that the HSM technique circumvents the notorious Gibbs phenomenon, which limits the accuracy that classic numerical methods reach.
% Consequently, the HSMs deliver a byproduct parameterization of the discontinuity manifold.
% 
Numerical experiments, addressing linear and non-linearly propagating shocks, demonstrate the strong approximation power of the HSM technique.
% In this paper, we introduce a novel approach that uses \emph{Hybrid Surrogate Models} (HSMs) for efficiently solving partial differential equations with non-smooth solutions, in particular solutions with finitely many jump discontinuities. The method relies on posing the PDE problem in a variational sense and using the HSMs as a surrogate for approximating the solution. We demonstrate with numerical experiments that the HSMs mitigate the Gibbs phenomenon while having a good approximation power, furthermore, it gives us as a subproduct a parametrization for the discontinuity manifold.
}


\section{Introduction}\label{sec1}
Solving partial differential equations (PDEs) with discontinuous solutions is crucial in various fields, including droplet physics \cite{WANG2022102684}, shock formation in hyperbolic systems \cite{Nakane+1994+199+204}, and singularity formation in magneto-hydrodynamics \cite{CHAE20161009}.
The challenge lies in the inherent difficulty of numerically approximating non-continuous functions, as smooth approximations struggle to resolve discontinuities, due to the notorious Gibbs phenomenon \cite{jerri1998gibbs}. Gibbs phenomenon manifests as spurious oscillatory behavior in various approximations, such as Fourier series \cite{zygmund_2003}, polynomials \cite{approx_th}, and even neural networks \cite{NN_ref}.
%crack propagation in mechanics \cite{Rege_2017}
Classic numerical methods, including finite differences \cite{LeVeque2007FiniteDM}, finite elements  \cite{ern2004theory}, finite volumes \cite{eymard2000finite}, and particle methods \cite{li2007meshfree}, incorporate smooth approximations within their formulations, thereby inheriting the limitations associated with the Gibbs phenomenon. State-of-the-art methods, such as Discontinuous Galerkin methods (DGM) \cite{cockburn:hal-01352444} with non-linear limiters or unfitted DGM \cite{DGM}, have been proven to be effective in controlling spurious oscillations associated with high-order finite element schemes.
Recent advancements in machine learning demonstrate promising results by leveraging the flexibility of neural networks (NNs), approximating discontinuous PDE solutions \cite{tseng2023discontinuity}. 

Complementary to the aforementioned methods, here, we present the \emph{Hybrid Surrogate Models} (HSMs) as a novel approach, circumventing the Gibbs phenomenon. The HSMs inherently incorporate discontinuities into the surrogate model, combining the Heaviside function with multivariate polynomials. Remarkably, this eliminates the need for limiters or adaptive discretizations for suppressing the Gibbs phenomenon. 

This work is a follow-up of the \emph{Sobolev PINNs} (SC-PINNs) \cite{SC_PINN} and the \emph{Polynomial Surrogate Models} (PSMs) \cite{cardona2023learning}, exploiting the therein delivered technique of \emph{Sobolev cubatures}.
\subsection{Basic Definitions and Notation}\label{sec:not}
We denote with $\mathcal{B}=\Omega\times (0,T)$ the product space of the $d$-dimensional open hypercube  $\Omega = (-1,1)^d\subseteq \mathbb{R}^d$, $d\in \mathbb{N}$ and a 
temporal domain $[0,T]\subseteq \mathbb{R}_{\geq0}$, $T>0$, and with  $\overline{\mathcal{B}}$ its closure. 
lWe assume the discontinuities $S_i \subseteq \R^{d+1}$ to be $d$-dimensional regular (smooth)  manifolds given as hypergraphs  
$S_i = \{(x,t)\in \overline{\mathcal{B}}: t = s_i(x)\}$, with $s_i: \overline{\Omega} \to [0,T]$ to be at least of class $s \in C^1(\overline{\Omega}, [0,T])$, such that $L_i(x,t) = s_i(x) -t$ is a regular level-set function, i.e., $\nabla_{x,t} L(x,t) \neq 0$, for all $(x,t) \in L_i^{-1}(0) = S_i$, $i=1,\dots, l$.

Let $BV(\mathcal{B})$ be the space of \emph{Bounded Variations} \cite{Kanwal2004}. For the $i$-th parametrization $s_i$ we associate the Heaviside function $H_{s_i}\in BV(\mathcal{B})$ as: 
\begin{equation}\label{eq:HS}
    H_{s_i}(x,t)= \begin{cases}
        1 &\,,  t\leq s_i(x)\\
         0 &\,, t> s_i(x)
         \end{cases} \,.
\end{equation}
We set $\mathbf{S}=\bigcup\limits_{i=1}^lS_i$ and  $\mathfrak{h}^TH_\mathbb{S}=\sum\limits_{i=1}^lh_iH_{s_i}$, where $\mathfrak{h} = \{h_i\}_{i=1}^l \subseteq \R$ denotes the constant heights. We define the space of functions with $l$ jump discontinuities at $\mathbb{S}$ as 
\begin{equation}
M_l(\mathcal{B})=\{f\in BV(\mathcal{B}): f=g +\mathfrak{h}^TH_\mathbb{S}, \textrm{ for some }g\in C^0(\mathcal{B})\}\,.
\end{equation} 

We denote with $\Pi_{n,d}$ the space of $d$-variate polynomials of maximum $l_\infty$-degree $n \in \N$, spanned by some polynomial basis $\{q_\alpha\}_{\alpha \in  A_{n,d}}$. Hereby,  $A_{n,d}=\{\alpha \in \mathbb{N}^d : \|\alpha\|_\infty \leq n\}$ denotes the associated multi-index set, which we assume to be lexicographically ordered.  In particular, we consider three possible bases \cite{approx_th}:
\begin{enumerate}[label=\roman*)]
    \item The canonical polynomials $x^\alpha= x^{\alpha_1}\cdots,x^{\alpha_d}$, $\alpha \in A_{n,d}$.
    \item The Chebyshev polynomials of the first kind 
$T_\alpha = T_{\alpha_1}\cdots, T_{\alpha_d}$, $\alpha \in A_{n,d}$,  with $T_k(x) = \mathrm{arccos}(k\cos(\pi x))$, $x \in [-1,1]$, for $k \in \N$.
\item The Lagrange basis $L_{\alpha}$ with $L_{\alpha}(p_\beta) = \delta_{\alpha,\beta}$, defined with respect to a given set of interpolation nodes $p_\beta \in P_{n,d}\subseteq[-1,1]^d$, both enumerated by $\alpha,\beta \in A_{n,d}$.
\end{enumerate}
\subsubsection{Sobolev Cubatures}\label{sec:SOB}
We shortly recap the notion of \emph{Sobolev Cubatures} from our prior work 
 \cite{SC_PINN,cardona2023learning}:  Let $P_{n,d+1}= P_{n,d}\times P_T\subseteq \mathcal{B}$, where $P_T\subseteq [0,T]$ denotes the re-scaled 1D Legendre grid. 
We denote with $\{L_\alpha\}_{\alpha\in A_{n,d+1}}\in \Pi_{n,d+1}$ the Lagrange basis with respect to $P_{n,d+1}$ and with $\omega_\alpha = ||L_\alpha||^2_{L^2(\mathcal{B})}$ the associated Gauss-Legendre weights. We introduce  $\mathbb{W}_{n,d+1}= \textrm{diag}(\{\omega_{\alpha}\}_{\alpha\in A_{n,d+1}}) \in \R^{|A_{n,d+1}|\times|A_{n,d+1}|}$ as the weight matrix and $\mathbb{D}_x,\mathbb{D}_{t} \in \mathbb{R}^{|A_{n,d+1}|\times |A_{n,d+1}|}$ as the polynomial differential matrices with respect to the Lagrange basis. For $\beta \in \N^{d+1}$, with $l_1$-norm $\|\beta\|_1 = k \in \N$ we introduce the $\beta$-th differential operator and the $k$-th weight matrix by 
\begin{equation}
    \mathbb{D}^\beta=\prod\limits_{i=1}^{d+1} \mathbb{D}_{\beta_i}\,, \quad \mathbb{W}^k_{n,d+1}=\sum\limits_{||\beta||_{1}\leq k}(\mathbb{D}^\beta)^T\mathbb{W}_{n,d+1}(\mathbb{D}^\beta)\,.
\end{equation}The \emph{Sobolev cubatures}, approximating the $H^k$-norm $\|\cdot\|_{H^k(\mathcal{B})}$ result as:
\begin{equation}\label{eq:sob_cub}
    ||u||_{H^k}(\mathcal{B}) \approx \mathfrak{u}^T\mathbb{W}^k_{n,d+1}\mathfrak{u}\,, \quad u \in H^k(\mathcal{B})\,,
\end{equation}
where $\mathfrak{u}=u(P_{n,d+1})\in \mathbb{R}^{|A_{n,d+1}|}$ denotes the evaluation vector.

\subsection{Contribution}
We introduce the \emph{Hybrid Surrogate Models} (PSMs) $\psi_{\theta}\in \Psi_{\theta}(\mathcal{B})$, parametrized by coefficients $\theta \in \Theta$, approximating functions with $l \in \N$ jump discontinuities. The HSMs are given by combining a smooth surrogate $Q_{\xi} \in \Pi_{n,d+1}$ with Heaviside functions \cite{Kanwal2004} $H_s\in BV(\mathcal{B})$. The HSMs induce the construction of a "stitching" map $\Phi : M_l(\mathcal{B})\rightarrow C^0(\mathcal{B})$ that subtracts the detected byproduct-parametrization of the discontinuities from the solution, resulting in a continuous output that can be approximated with smooth polynomial surrogates. 

The HSM solutions are computed by finding the optimal coefficients $\theta\in\Theta$, minimizing an associated discretized loss $\mathcal{L}^n:\Theta \rightarrow \R_{\geq 0}$, formulated for detecting the discontinuities. Here, we focus on two challenges: 

\begin{itemize}
    \item \textbf{Reconstruction problem:} Given a function, $g\in M_l(\mathcal{B})$ we seek for a  HSM  approximate $g \approx \psi_{\theta}\in  \Psi_\theta(\mathcal{B}) \subseteq M_l(\mathcal{B})$, minimizing an associated discretized reconstruction loss $\mathcal{L}_0^n:\Theta\rightarrow\R_{\geq0}$.
    \item \textbf{1D transport equation problem:} Given an initial condition, $u_0\in M_l(\Omega)$ we want to find an HSM approximate $\psi_{\theta}\in \Psi_\theta(\mathcal{B})$ minimizing a discretized PDE loss $\mathcal{L}^n_{\textrm{PDE}}:\Theta\rightarrow \R_{\geq0}$ associated to the transport equation with Dirichlet boundary conditions $g\in C^0(\partial\Omega\times (0,T])$ and constant velocity $v\in \R$, given by:
    \begin{equation}\label{eq:tr_pde}
          \partial_t u +v\partial_x u = 0\,, \quad \textrm{in}\quad \mathcal{B}.
    \end{equation}
\end{itemize}
The capability of the HSM technique is demonstrated in numerical experiments.
\section{Hybrid Surrogate Models HSMs}\label{sec:HSM_app}
We detail the notion of the \emph{Hybrid Surrogate Models} and introduce the notion of distributional derivative for the Heaviside function: To start with we consider 
$$Q_{\xi}(x,t)= \sum\limits_{\alpha \in A_{n,d+1}}\xi_\alpha T_\alpha(x,t) \in \Pi_{n,d+1}\,, \quad  s_i(x)= \sum\limits_{\alpha\in A_{m,d}}(C_i)_\alpha x^\alpha\in  \Pi_{m,d}$$ 
a polynomial $Q_\xi$ in \emph{Chebyshev} basis and the approximation $s$ of the parametrizations  $s_i^*\in C^1(\overline{\Omega};[0,T])$ in canonical basis.  We define the space of HSMs $\Psi_\theta(\mathcal{B})\subseteq M_l(\mathcal{B})$ as:
\begin{equation}\label{def:HSM_space}
    \Psi_\theta(\mathcal{B})=\left \{ Q_{\xi} + \mathfrak{h}^TH_{\mathbb{S}} : Q_{\xi} \in \Pi_{n,d+1}\right\}\,,
\end{equation}
where $\mathfrak{h}^T H_{\mathbb{S}}
= \sum\limits_{i=0}^lh_iH_{s_i} \in BV(\mathcal{B})$ with constant hights $\mathfrak{h}= \{h_i\}_{i=1}^l\subseteq \R$. We assume the 
elements $\psi_{\theta}\in   \Psi_\theta(\mathcal{B})$ to be parametrized by $\theta=\{\xi, \beta, C\}\in \Theta \subseteq \R^p$, with $p= |A_{n,d}|+l+l\cdot|A_{m,d}|$.

\subsection{Distributional Derivative} 
For leveraging the HSMs solving PDEs, we define the differential operators acting on the HSMs in the distributional sense. When restricting to the 1D setting $d =1$ it suffices to provide the notion only with respect to time $t\in [0,T]$. The definition involves 
the spaces of test functions $\mathcal{D}(\mathcal{B})$ and the space of distributions $\mathcal{D}'(\mathcal{B})$ \cite{Kanwal2004}.
\begin{lemma}\label{lemm:dist_der}
    Let $d=1$, $\Omega = (-1,1)$ and $H_s\in BV(\mathcal{B})$ be the Heaviside function with jump discontinuity parametrized by  $s\in C^1(\overline{\Omega};[0,T])$, Eq.~\eqref{eq:HS}. 
    \begin{enumerate}[label=(\roman*)]
        \item\label{L1}  Let 
    $F\in \mathcal{D}'((0,T))$  be the distribution
on time dependent functions, given by   $$F[\phi(x,\cdot)] = \int\limits_{0}^{T}H_s(x,t)\phi(x,t)dt\,,  \quad \phi(x,\cdot) \in \mathcal{D}((0,T))\,, \quad x \in \Omega\,.$$ 
    Then the distributional derivative $F'[\phi(x,\cdot)] = - F[\partial_t\phi(x,\cdot)]$ is equivalently given by 
    \begin{equation*}
        F'[\phi(x,\cdot)] = \phi(x,s(x)) = \int\limits_{0}^T \phi(x,t)\delta_s(t)dt\,,
    \end{equation*}
    implying the distributional derivative $\partial_tH_s(x,t) =\delta_s(t)$ to be given by the  Dirac function.
    \item\label{L2} The distributional derivatives $\partial_xH_s(x,t), \partial_tH_s(x,t)$ satisfy the relation
    \begin{equation}
        \partial_x H_s(x,t) = \partial_tH_s(x,t)\partial_x s(x)\,,
    \end{equation}
    \end{enumerate}
\end{lemma}
\begin{proof} \ref{L1} is a standard result and can be found in \cite{Kanwal2004}. To show \ref{L2} we 
    note that by our assumptions on the parametrization $s\in C^1(\overline{\Omega};[0,T])$, section \ref{sec:not}, the inverse function theorem \cite{bams/1183549049} applies. Hence,  
     $x= s^{-1}(t)$ for all $(x,t)\in S$. Consequently, when applying Fubini's theorem and \ref{L1}, the following holds:
    \begin{align}
    \int\limits_{\Omega}\int\limits_{0}^T\partial_xH_s(x,t)\phi(x,t)dtdx =  -\int\limits_{0}^T\phi(s^{-1}(t),t)dt = -\int\limits_\Omega\phi(x,s(x))\partial_xs(x)dx\,.
    \end{align}
    On the other side, by \ref{L1} the distributional derivative  
    $\partial_t H_s(x,t)$ is given by
    \begin{align}
    \int\limits_{\Omega}\int\limits_{0}^T\partial_tH_s(x,t)\phi(x,t)dt =  -\int\limits_\Omega \phi(x,s(x))dt,
    \end{align}
    yielding the claim.
\end{proof}

We provided all the ingredients for addressing the reconstruction problem. 
\section{Reconstruction Problem and Discretized Loss}
We state a variational optimization problem, providing HSM approximates of non-continuous functions $g\in M_l(\mathcal{B})$. 
To do so, we introduce reconstruction loss functional $\mathcal{L}_0:\Psi_\theta(\mathcal{B})\rightarrow\mathbb{R}_{\geq 0}$, by
\begin{equation}\label{eq:rec_var}
    \inf\limits_{\psi_\theta \in \Psi_\theta(\mathcal{B})}\mathcal{L}_0[\psi] \textrm{, with }\mathcal{L}_0[\psi_\theta]= \int\limits_{\Omega}\int\limits_0^T \biggr(\psi_\theta(x,t) - g(x,t)\biggr)^2 dtdx.
\end{equation}
For the sake of simplicity, we assume that $g\in M_l(\mathcal{B})$ possesses  only one single discontinuity, $l=1$, parametrized by $s^*\in C^1(\overline{\Omega};[0,T])$ with constant height $h^*\in \R$. Given $\psi_\theta\in \Psi(\mathcal{B})$ and its approximate $s \approx s^*$, $s \in \Pi_{m,d+1}$ we define the subdomains 
\begin{equation}
    U_{-}(s)= \{(x,t)\in \mathcal{B}: t\leq  s(x)\},\,\,\,\, 
    U_{+}(s)= \{(x,t)\in \mathcal{B}: t > s(x)\}
\end{equation} 
Applying the definition of $\Psi_\theta(\mathcal{B})$ from Eq.~\eqref{def:HSM_space}, the reconstruction loss $\mathcal{L}_0$ is equivalent to: 
\begin{equation}\label{eq:seg_loss}
    \mathcal{L}_0[\psi_{\theta}] =  \int\limits_{U_{-}} (Q_{\xi}+h-g)^2dtdx + \int\limits_{U_{+}} (Q_{\xi}-g)^2dtdx.
\end{equation}
 We exploit the Sobolev cubatures for discretizing the loss.
 \begin{defi}[discretized reconstruction loss] Let $\mathfrak{u}_{\pm}(s)=\{(x_\alpha, t_\alpha)\in P_{n,d+1}: (x_\alpha,t_\alpha)\in U_{\pm}(s)\}$ be the left/right training sets, 
 $A^{\pm}_{n,d+1}=\{\alpha \in A_{n,d+1}: (x_\alpha,t_\alpha)\in U_{\pm}(s)\}$, and   
 $\mathbb{T}^{\pm}_{n,d+1}(s)=(T_i(\mathfrak{u}_{\pm}(s)))_{i\in {A_{n,d+1}}}$ denote the left/right Vandermonde matrices with respect to the \emph{Chebyshev} basis,  and $\mathbb{W}^{k\pm}_{n,d+1}=(\mathbb{W}^k_{n,d+1})_{i,j\in A^\pm_{n,d+1}}$ be the corresponding left/right Sobolev weight matrices.
The discretized reconstruction loss $\mathcal{L}^{n}_0:\Theta\rightarrow\mathbb{R}_{\geq 0}$, 
Eq.~\eqref{eq:seg_loss} with Sobolev cubatures of order $k>0$ is defined as 
%\begin{lemma}\label{eq:schk_disc}
%    \MH{"Using the left and right definitions, the discrete reconstruction loss is given by".. this is no valid sentence / formulation.}
\begin{equation}\label{lemm:disc_vp}
    \mathcal{L}^n_0[\theta] = (\mathfrak{r}^{-}_{n,d+1})^T\mathbb{W}^{k-}_{n,d+1}\mathfrak{r}^{-}_{n,d+1} + (\mathfrak{r}^{+}_{n,d+1})^T\mathbb{W}^{k+}_{n,d+1}\mathfrak{r}^{+}_{n,d+1}\,,     
\end{equation}
where the left and right residuals are given by $\mathfrak{r}^{-}_{n,d+1}= \mathbb{T}^{-}_{n,d+1}\xi +h - g(P_{n,d+1}^{-})$, $\mathfrak{r}^{+}_{n,d+1}= \mathbb{T}^{+}_{n,d+1}\xi - g(P_{n,d+1}^{+})$.
     \label{def1}
 \end{defi}
 We further extend the concept by addressing the PDE problem below. 

\section{Transport Equation in 1D with Discontinuous Solution}
We focus on exploiting the HSMs for solving the 1D transport equation, Eq.~\eqref{eq:tr_pde}, with an initial condition $u_0\in M_l(\Omega)$, possessing one single discontinuity at  $x_0\in \Omega$, constant velocity $v\in \R$, and Dirichlet boundary conditions $g\in C^0(\partial\Omega\times[0,T])$. The corresponding solution $u \in M_l(\mathcal{B})$ has a single jump discontinuity, $l=1$, parametrized by $s^*\in C^1(\overline{\Omega};[0,T])$, and a constant height $h^*\in \R$.

Stating the problem in a distributional sense is to seek for a solution $u\in M_l(\mathcal{B})$, minimizing the loss
\begin{equation}\label{eq:dist_pde}
    \mathcal{L}_\textrm{PDE}[u]=\sup\limits_{\phi \in D(\mathcal{B})}\biggr(\int\limits_{\Omega}\int\limits_{0}^T(\partial_t u (x,t)+v\partial_x u(x,t))\phi(x,t) dt dx\biggr)^2.
\end{equation}
%Hereby, we assume a constant velocity $v\in\mathbb{R}$, a non-continuous initial condition $u_0\in BV(\Omega)$ and Dirichlet boundary conditions $g\in C^0(\partial\Omega\times(0,T])$. 
%\MH{ State the equation explicitly !}
We propose to derive an  approximate surrogate solution $\psi_\theta\in \Psi_\theta(\mathcal{B})$ by minimizing the discrete loss $\mathcal{L}^n:\Theta\rightarrow \R_{\geq 0}$, given by
\begin{equation*}
    \mathcal{L}^n[\theta]= \mathcal{L}^n_\textrm{PDE}[\theta] + \mathcal{L}^n_{\partial}[\theta]+ \mathcal{L}^n_{0}[\theta]\,,
\end{equation*}
where $\mathcal{L}^n_{\partial}[\theta]$, $\mathcal{L}^n_{0}[\theta]$ are the Sobolev cubature discretizations of the boundary and reconstruction losses
\begin{equation*}
\mathcal{L}_{\partial}[\psi_\theta]= \int\limits_{\partial\Omega}\int\limits_{0}^T(\psi_\theta - u_0)^2dtdx\,, \quad 
\mathcal{L}_{0}[\psi_\theta]= \int\limits_{\Omega}(\psi_\theta(x,0) - g(x))^2dx\,.
    %&\mathcal{L}_\partial[\theta]= (\mathfrak{r}_{n,d})^T\mathbb{W}^{-}_{n,d}\mathfrak{r}^{-}_{n,d} + (\mathfrak{r}^{+}_{n,d})^T\mathbb{W}^{+}_{n,d}\mathfrak{r}^{+}_{n,d}\\&
    % \mathcal{L}_0[\theta]= (\mathfrak{s}_{n,d})^T\mathbb{W}_{n,d}\mathfrak{s}_{n,d}
\end{equation*}
The discrete PDE loss, Eq.~\eqref{eq:dist_pde}, results when restricting the space of test functions to the polynomial space $\Pi_{n,d+1}$:
\begin{equation}\label{eq:dist_loss1}
    \mathcal{L}^n_{\mathrm{PDE}}[\theta]=\sum\limits_{\alpha \in A_{n,d+1}}\biggr(\int\limits_{\Omega}\int\limits_{0}^T(\partial_t \psi_\theta +v\partial_x \psi_\theta)L_\alpha dt dx \biggr)^2 \,.
\end{equation}
The following Lemma, gives us an equivalent equation for the discrete PDE loss, by exploiting the structure of the HSM surrogates.
\begin{lemma}\label{thm:dist}
     Let $\psi_\theta\in\Psi_\theta(\mathcal{B})$ parametrized by $\theta=\{\xi,h,C\}$, then the discretized PDE loss associated to Eq. \eqref{eq:dist_loss1} is equivalent to: 
    \begin{align}\label{eq:sd_loss}
        \mathcal{L}^n_{\mathrm{PDE}}[\theta] =& \sum\limits_{\alpha \in A_{n,d+1}}\biggr(\int\limits_{\Omega}\int\limits_{0}^T(\partial_tQ_{\xi} + v\partial_xQ_{\xi})L_\alpha(x,t)dtdx\\& -h\int_\Omega((1+v\partial_x s(x))L_\alpha(x,s(x)))dx \biggr)^2 \label{eq:DP}
    \end{align}
\end{lemma}

\begin{proof}
 Due to the linearity of the PDE operator $\partial\cdot +v\partial_x\cdot$, the definition of the HSMs, and the smoothness of $Q_{\xi}\in \Pi_{n,d+1}$ we obtain the first term, $\partial_tQ_{\xi} + v\partial_xQ_{\xi}$, in Eq.~\eqref{eq:sd_loss} in the strong sense. The discontinuous part, Eq.~\eqref{eq:DP} follows from Lemma~\eqref{lemm:dist_der}.
\end{proof}   

\begin{defi}[discretized PDE loss] Given the ingredients of Eq.~\eqref{eq:dist_pde} we define the discretized PDE loss, Eq. \eqref{eq:sd_loss}, in terms of the Sobolev cubatures of order $k=0$ as:
        \begin{align*}
    \mathcal{L}^n_{\textrm{PDE}}[\theta] = \sum\limits_{\alpha \in A_{n,d+1}}\biggr(&\sum_{\gamma \in A_{n,d+1}}(\mathbb{D}_t\mathbb{T}_{n,d+1}\xi + v\mathbb{D}_x\mathbb{T}_{n,d+1}\xi)_\gamma \omega_\gamma\\&
    -h\sum_{\gamma \in A_{n,d}}\bigr((1+v\partial_xs(P_{n,d}))L_\alpha(s(P_{n,d}))_\gamma\bigr)\omega_\gamma \biggr)^2\,,
\end{align*}
where $\mathbb{T}_{n,d+1}=(T_\alpha(P_{n,d+1}))_{\alpha\in A_{n,d}}\in \R^{|A_{n,d+1}|\times|A_{n,d+1}|}$ denotes the Vandermonde matrix with respect to the Chebyshev polynomials and the Legendre grid.
\label{def2}
\end{defi}
Definition \ref{def2} can be extended to Sobolev cubatures of arbitrary order $k\in \Z$ \cite{cardona2023learning}. The setup $k=-1$ is used in Experiment~\ref{Ex:2}, while Definition~\ref{def1} with $k=0$ is the setting of our implementation for Experiment~\ref{Ex:1}.
\section{Numerical Experiments}
We designed experiments, addressing the reconstruction problem of a non-linearly propagating shock, and the transport equation in 1D with discontinuous initial conditions. 

We compare the HSMs with the methods mentioned in the introduction: PSMs \cite{cardona2023learning}, SC-PINNs \cite{SC_PINN}, and standard NNs with ReLU and Sine activation functions. All models were trained over the same number of points. For HSMs, SC-PINNs, and PSMs we used Legendre grids of degree $n_x,n_t = 72$ in the $x$ and $t$ direction, respectively, while for the standard NNs, the training points were randomly sampled from $\mathcal{B}$. The NN-based models have an architecture of $5$ fully connected layers with $50$ activation functions each. The errors are measured in the $l^1$-norm for 
$N=300^2$ equidistant points $\{x_i,t_i\}_{i=1}^N\subseteq \mathcal{B}$.
 
% $$\epsilon_1=\frac{1}{N}\sum\limits_{i=1}^N|\hat{u}_\theta(x_i,t_i) - g(x_i,t_i)|\,, \,\, \text{for} \,\,  N=300^2 \,\, \text{equidistant points} \,\,\  \{x_i,t_i\}_{i=1}^N\subseteq \mathcal{B}\,.$$
\begin{figure}[h]
\hspace{10pt}\subfloat[Reconstruction Problem]{
\label{fig:gt_1a}
\centering
\includegraphics[width=0.33\textwidth]{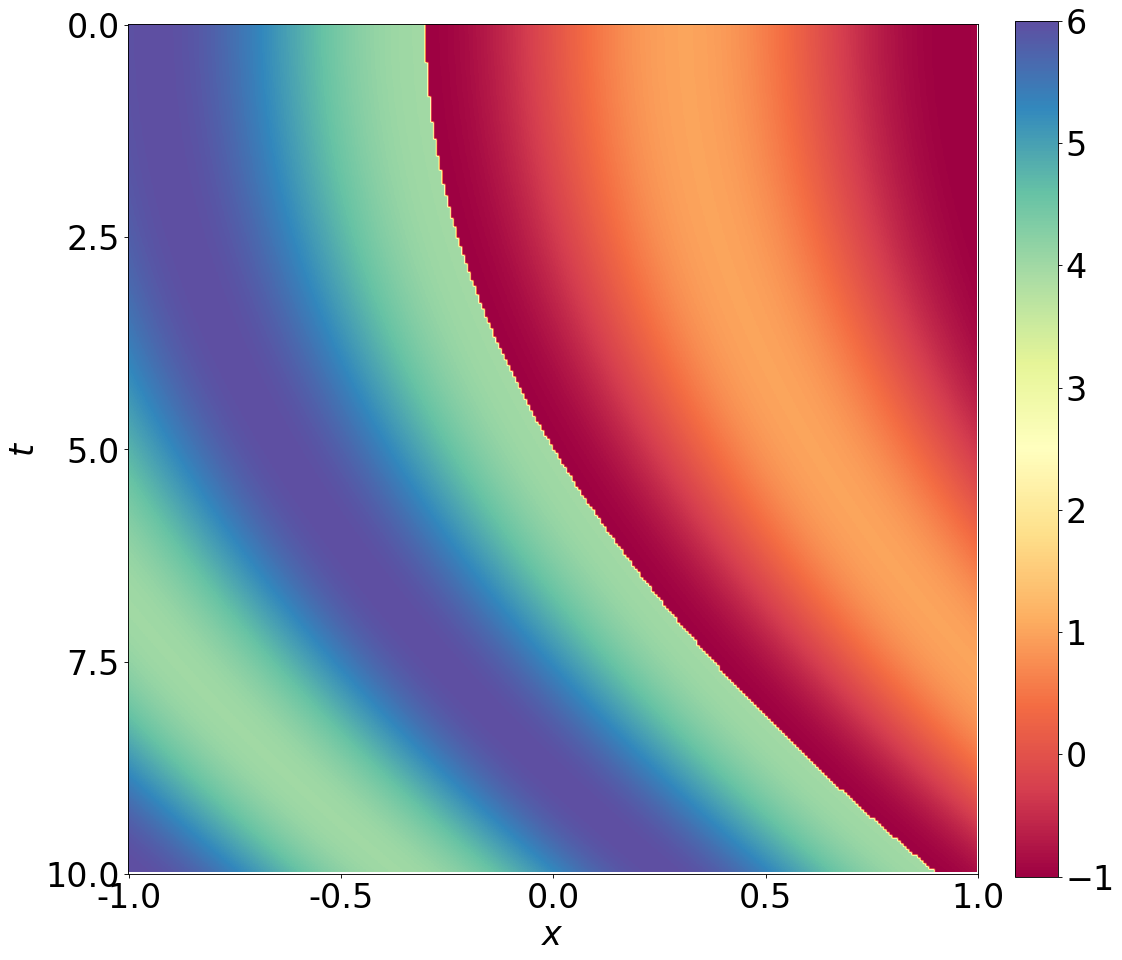}
}
%\begin{subfigure}[b]{.30\textwidth }%
\hspace{40pt}\subfloat[Transport Equation]{
\label{fig:gt_1b}
\centering
\includegraphics[width=0.33\textwidth]{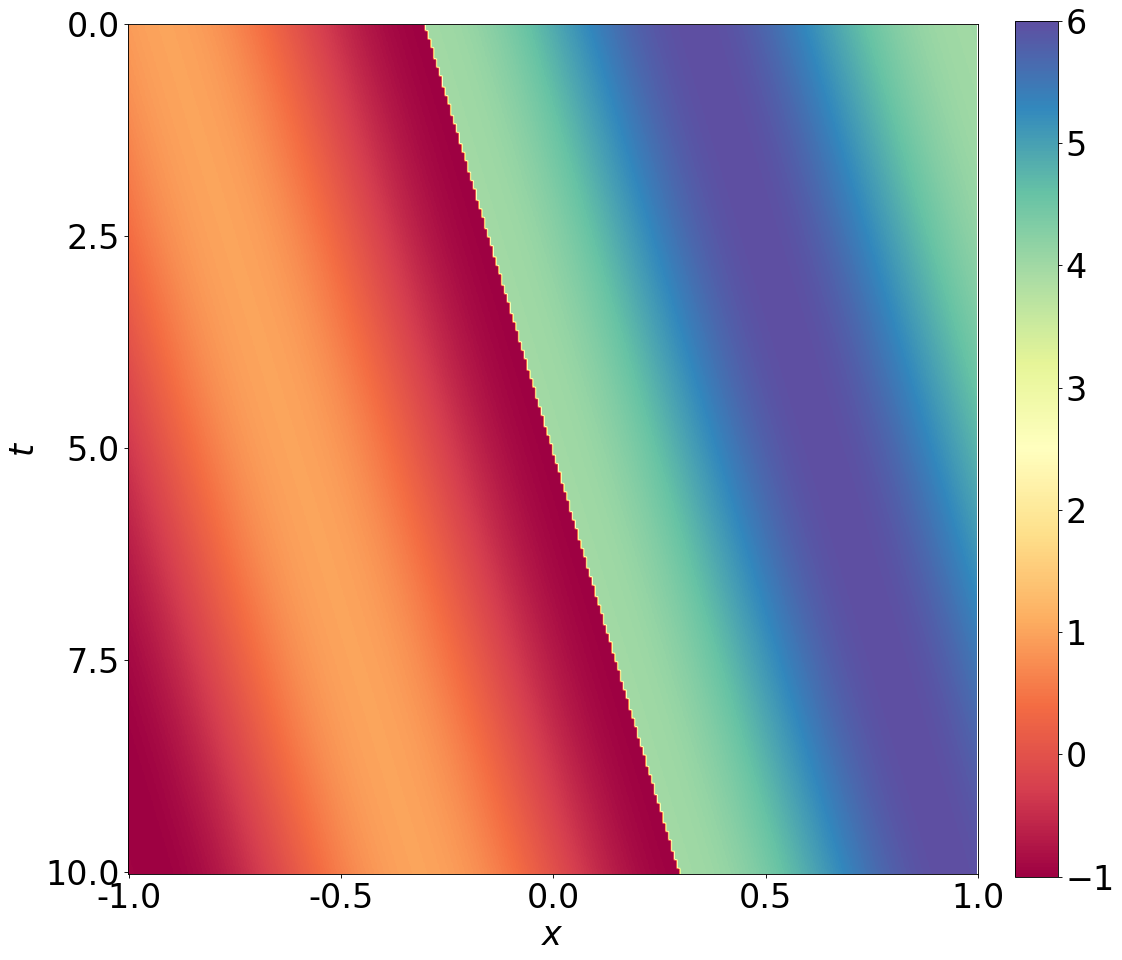}
}
\caption{Ground truth for the numerical experiments}
\end{figure}
\begin{experiment}\label{Ex:1}(\textbf{Reconstruction Problem of a Discontinuous Function})
 We consider the setting with $d=1$ of approximating a discontinuous function $g\in BV(\mathcal{B})$ with a single jump discontinuity, $l=1$, constant height $h^* = 5$, and analytical smooth part given by $\Phi(g)(x,t) = \sin(5x-f(t))$ with $f(t)=0.3t^2-0.3$.
 \end{experiment}

\begin{figure}[h]
\subfloat[Predictions at $t=0.33$]{
\label{fig:approx_1a}
\centering
\includegraphics[width=0.4\textwidth]{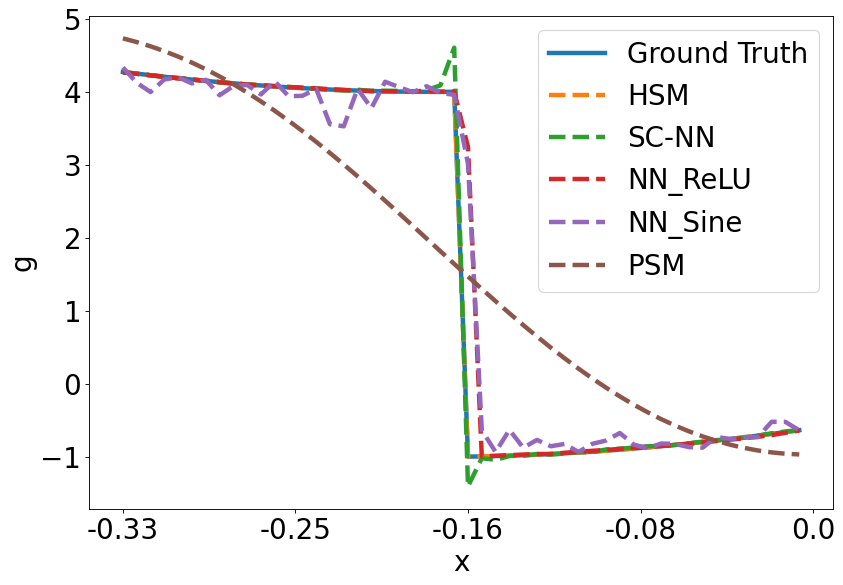}
}
\hspace{20pt}\subfloat[Approximation errors in $\mathcal{B}$]{
\label{fig:approx_1b}
\centering
\begin{tabular}{lll}
method  & $l_1$-error & runtime (s)\\
\hline
  PSM & $3.74 \cdot10^{-2}$ & $\bf t = 2.37$ \\
  SC-PINN  & $3.79 \cdot10^{-3}$ &  $t = 203.50$\\
  NN-ReLU & $4.92\cdot 10^{-3}$ & $t = 121.60$ \\
  NN-Sine & $2.83\cdot10^{-2}$ &  $t = 124.46$\\
  HSM & $\bf 7.33\cdot10^{-5}$ &  $t = 19.38$\\ 
   & &  \\ 
   & &  
\end{tabular}
}
\caption{Approximation errors for the reconstruction problem}\label{tab:rec_p}
\end{figure}

Table~\ref{fig:approx_1b} reports the errors. HSM reaches   
at least two orders of magnitude better accuracy, while requiring one order of magnitude less computational time compared to the NN based models. When considering the solutions at $t=0.33$ in 
Fig.~\ref{fig:approx_1a}, we observe that PSM, NN\_Sine and SC\_PINN run into Gibbs phenomenon. The NN\_ReLU avoids the Gibbs phenomenon but fails to capture the discontinuity accurately, unlike the HSM that fits the discontinuity optimally.

We continue our investigations:

\begin{experiment}\label{Ex:2}(\textbf{Transport equation in 1D with discontinuous initial condition}) We keep the design of Experiment \ref{Ex:1} and 
consider the 1D transport equation, with an initial condition possessing a discontinuity at $s_0=-0.3$, a smooth part given by $\Phi(u_0(x))= \sin(5x)$ and a constant height $h^*= 5$. We assume Dirichlet boundary conditions given by the analytical solution $u(x,t) = u_0(x-vt)$ and velocity $v= 1$. 
\end{experiment}

\begin{figure}[h!]
\subfloat[Predictions at $t=0$]{
\label{fig:pde_1a}
\centering
\includegraphics[width=0.4
\textwidth]{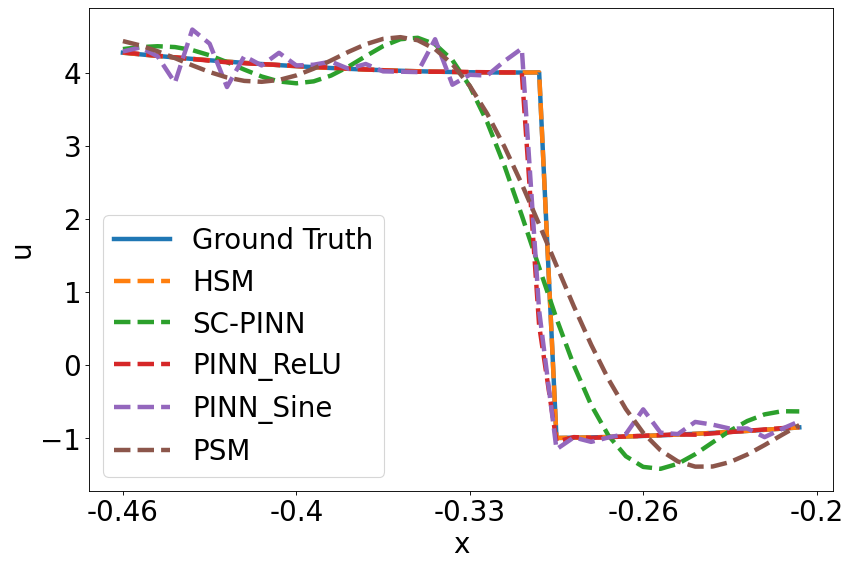}
}
\hspace{20pt}\subfloat[Approximation errors in $\mathcal{B}$]{
\label{fig:pde_1b}
\begin{tabular}[t]{lll}
method & $l_1$-error& runtime (s)\\
\hline
  PSM & $3.77 \cdot10^{-2}$ & $\bf t = 0.89$ \\
  SC-PINN  & $1.44 \cdot10^{-3}$ &  $t = 214.47$\\
  NN-ReLU & $6.07\cdot 10^{-3}$  & $t = 342.45$ \\
  NN-Sine & $1.75 \cdot10^{-2}$  &  $t = 404.78$\\
  HSM & $\bf 1.45 \cdot10^{-4}$ &  $t = 29.00$\\
& & \\
& &
\end{tabular}%
}
\caption{Approximation of the 1D transport equation}\label{NVE}
\end{figure} 
Table~\ref{fig:pde_1b} reports the errors. As expected all the methods except for the NN\_ReLU and the HSM run into Gibbs phenomenon. All the methods are outperformed by the HSM, both in accuracy and in runtime. Fig.~\ref{fig:pde_1a} clearly shows the HSM, circumventing the Gibbs phenomenon optimally. 

\section{Conclusions}
% In conclusion, this paper presents a novel approach for solving partial differential equations with non-smooth solutions characterized by finitely many jump discontinuities. 
The \emph{Hybrid Surrogate Models (HSMs)} showed remarkable efficiency and approximation power when dealing with discontinuous functions, avoiding the Gibbs phenomenon while providing an approximated parametrization of the discontinuity as a byproduct. 
% the dise challenging solutions by exploiting the approximation power of the polynomial surrogates and incorporating the structure of the discontinuity on the Heaviside function. By formulating the PDE problem in a variational framework and employing HSMs as surrogates, we have successfully addressed the Gibbs phenomenon while achieving up to two orders of magnitude improvement in the $l^1$-norm of the error compared to Neural Network-based methods or global polynomials. 
% 
Extending the present 1D implementations to multiple shocks and higher dimensions
% Furthermore, unlike the other methods considered, this approach provides a parametrization for the discontinuity manifold as a valuable byproduct. This innovative technique 
holds great promise for advancing the field of shock-type PDE solvers.
\begin{acknowledgement}
This research
was partially funded by the Center of Advanced Systems Understanding (CASUS), which is financed by Germany’s Federal Ministry of Education and Research (BMBF) and by the Saxon Ministry for Science, Culture, and Tourism (SMWK) with tax funds on the basis of the budget approved by the Saxon State Parliament.
\end{acknowledgement}
\bibliographystyle{abbrv}
\bibliography{REF}

\end{document}